\documentclass[12pt]{article}
\usepackage[centertags]{amsmath}
\usepackage{color,latexsym,amsfonts,amssymb}
\def\ignore#1{}

\newtheorem{thm}{Theorem}[section]

\newtheorem{cor}[thm]{Corollary}
\newtheorem{lem}[thm]{Lemma}
\newtheorem{Remark}[thm]{Remark}

\def \B {{\cal B}}  \def \F {{\cal F}} \def\ci{{\cal I}} \def \L {{\cal L}} \def\R{\mathbb{R}} 
\def \A {{\cal A}}
\def\Z{\mathbb{Z}}
\def \leo {\left(}
\def \rio {\right) }

\newcommand{\qed}{\hbox{\vrule width 5pt height 5pt depth 0pt}}
\makeatletter
\renewcommand\theequation{\thesection.\@arabic\c@equation}
\makeatother

\newcommand{\beas}{\begin{eqnarray*}}
\newcommand{\enas}{\end{eqnarray*}}

\newcommand{\bea}{\begin{eqnarray}}
\newcommand{\ena}{\end{eqnarray}}
\newcommand{\eq}{\begin{equation}}
\newcommand{\en}{\end{equation}}

\def\a{{\alpha}}
\def\b{{\beta}}

\def\e{\varepsilon}
\def\D{\Delta}

\def\nexto{\kern -0.54em}

\def\bone{{\bf 1}}
\def\bZ{{\bf Z}}

\def\l{{\lambda}}
\newcommand{\mean}{{\rm {I\ \nexto E}}}
\def\ex{\mean}

\newcommand{\Pn}{{\rm Pn\,}}

\newcommand{\PBD}{{\rm PBD}}
\def\prob{{\rm {I\ \nexto P}}}

\def\Ref#1{(\ref{#1})}

\newcommand{\tg}{{\tilde g}}

\begin{document}

\title{Polynomial birth-death distribution approximation in Wasserstein distance}

\author{ Aihua Xia\footnote{E-mail: xia@ms.unimelb.edu.au}
 \ and \  Fuxi Zhang\footnote{E-mail: zhangfxi@pku.edu.cn}
\\
Department of Mathematics and Statistics\\
The University of Melbourne\\
VIC 3010, Australia}
\date{{{21 November, 2008}\\
(First version: 24 November, 2006)}}

\maketitle

\begin{abstract}
The polynomial birth-death distribution (abbr. as PBD) on $\ci=\{0,1,2, ...\}$ or $\ci=\{0,1,2, ... ,m\}$ for some finite $m$ introduced in Brown \& Xia~(2001) is the equilibrium distribution of the birth-death process with birth rates $\{\alpha_i\}$ and death rates $\{\beta_i\}$,
where $\a_i\ge0$ and $\b_i\ge0$ are polynomial functions of $i\in\ci$. The family includes Poisson, negative binomial, binomial and hypergeometric distributions. In this paper, we give probabilistic proofs of various Stein's
factors for the PBD approximation with $\a_i=a$ and $\b_i=i+bi(i-1)$ in terms of the Wasserstein
distance. The paper complements the work of Brown \& Xia~(2001) and generalizes the work of Barbour \& Xia~(2006) where Poisson approximation ($b=0$) in the Wasserstein distance is investigated. As an application, we establish an upper bound for the Wasserstein distance between the PBD  and Poisson binomial distribution and show that the PBD approximation to the Poisson binomial distribution is much more precise than the approximation by the Poisson or shifted Poisson distributions.



\vskip12pt \noindent\textit{Running title\/}:
PBD approximation in Wasserstein distance
\end{abstract}

\section{Introduction and the main results.}
\setcounter{equation}{0}

{The cornerstone of the `law of small numbers' is the Poisson limit theorem which says that the total number of independent (or weakly dependent) rare events follows approximately the Poisson distribution. Prohorov~(1953), Hodges \& Le Cam~(1960) [also known as the Hodges-Le Cam theorem] and Chen~(1975)  can be regarded as the milestones for quantifying the Poisson limit theorem. However, Barbour and Hall~(1984) proved that the accuracy of Poisson approximation in total variation is determined by the rarity of the events and its order does not improve when the sample size increases. Inspired by the pioneering work in Poisson approximation, various attempts have been made to improve on the precision of approximation, leading to new forms of the `laws of small numbers', e.g., binomial [Ehm~(1991)], compound Poisson [Barbour, Chen \& Loh~(1992)], signed compound Poisson [Presman~(1983), Kruopis~(1986), \v Cekanavi\v cius~(1997)] and polynomial birth-death (abbr. as PBD) [Brown \& Xia~(2001)]. In particular, PBD approximation to the number of independent (or weakly dependent) rare events in total variation distance has been proved extremely accurate in Brown \& Xia~(2001), i.e. the errors of approximation decrease when the events become rarer and/or the sample size increases. On the other hand, all these approximations can also be considered in the context of the Wasserstein distance $d_W$ [see Barbour, Holst \& Janson~(1992), p.~13]: for any two distributions ${\cal P}$ and ${\cal Q}$ on $\Z_+:=\{0,1,2, ...\}$,
 $$
d_W ({\cal P}, {\cal Q}) = \sup_{f \in \F} \left| \int f d {\cal P}  - \int f d
{\cal Q} \right|,
 $$
where $\F$ is the set of functions $f : \Z_+ \rightarrow \R$, $|f(x) - f(y)|
\le |x-y|$ for $x,y\in\Z_+$. The metric $d_W$ characterizes the weak convergence and the convergence of the first absolute moments \ignore{and is a metric needed in the Hugarian construction in the studies of empirical processes} [see Shorack \& Wellner (1986), pp. 64--65].
Poisson approximation with respect to the Wasserstein distance is well documented in Barbour, Holst \& Janson~(1992), pp. 13--17 and Barbour \& Xia~(2006). In this note, we establish Stein's factors for PBD approximation in terms of the Wasserstein distance and the work generalizes the results in Barbour \& Xia~(2006).  To demonstrate the significance of our results, we prove an upper bound for the Wasserstein distance between the PBD  and Poisson binomial distribution. The bound is in the fashion of the Hodges-Le Cam theorem and it implies that the PBD approximation to the Poisson binomial distribution is much more precise than the approximation by the Poisson or shifted Poisson distributions.}

The {PBD} distribution introduced in Brown \& Xia~(2001) is the discrete distribution on $\ci=\Z_+$ or $\ci=\{0,1,2, ... ,m\}$ for some finite $m$ with probability function
\eq
\pi_n = \frac{\prod_{i=0}^{n-1} \a_i}{\prod_{i=1}^n \b_i} \left( 1 +
\sum_{j=1}^\infty \frac{ \prod_{i=0}^{j-1} \a_i}{\prod_{i=1}^j \b_i} \right)^{-1}, \
n\in\ci,
\label{equilibriumdis}
\en
where $\a_i\ge0$ and $\b_i{\ge}0$ are polynomial functions of $i\in\ci$. The framework for PBD approximation unifies
Poisson approximation [Chen~(1975)] adapted from Stein~(1972)
(in this case it is often appropriately called the Stein-Chen method), binomial approximation
[Ehm~(1991)], negative binomial approximation [Brown and Philips~(1999)] and some compound Poisson approximation [Barbour, Chen \& Loh~(1992)]. 
In fact, it is a routine exercise to check that a random variable $X$ follows the
distribution $\pi$ iff  
$$\mean[ \a_X g(X+1) - \b_X g(X)] = 0
 $$
 for every function $g$ satisfying $\mean[\b_X |g(X)|]<\infty$.
For any function
$f$ on $\ci$ satisfying $\sum_{i\in\ci} |f(i)|\pi_i<\infty$, one can recursively solve for the function $g_f$ such that
 \eq \label{Steinequ}
{\cal B} g_f(i) : = \a_i g_f(i+1) - \b_i g_f(i) = f(i) - \pi (f),
 \en
 where $\pi(f)=\sum_{i=0}^\infty f(i)\pi_i$. Then for any random variable $W$,
 \eq
| \mean {\cal B} g_f(W) | = | \mean f(W) - \pi (f) |
\label{AX01}
 \en
gauges the difference between $\L W${, i.e. the distribution of $W$,} and $\pi$ in terms of the test function $f$. To
estimate the {Wasserstein} distance between $\L W$ and $\pi$, it is sufficient to calculate the
supremum of the left hand side of \Ref{AX01} over the {test functions $f\in\F$}.

As in Brown \& Xia~(2001), we use $\PBD(\a;0,\b,1)$ to stand for the distribution $\pi$ when $\a_i = \a$ and $\b_i = \b i + i (i-1)$ with
$\a,\ \b\ge 0$. {Noting that the value of $g_f(0)$ has no contribution to the equation \Ref{Steinequ}, we take $g_f(0)=g_f(1)$ in this paper.} In estimating the
difference between the distribution of a random variable $W$ and $\PBD(\a;0,\b,1)$,
one often needs to resort to the first and second order Palm distributions of $W$ [Kallenberg~(1983), p.~103]. For the
distribution of the sum of independent Bernoulli variables, it is shown in Brown \&
Xia~(2001) that $\PBD(\a;0,\b,1)$ offers very good quality of approximation. In
fact, 
for random variables with well-behaved first and second order Palm distributions,
$\PBD(\a;0,\b,1)$ is a natural choice amongst many approximating distributions. 

Let $\|\cdot\|$ denote the supremum norm of the function $\cdot$ over its range.
Depending on the choices of values of $\a$ and $\b$, we can often show that 
$$
  |\mean\{\a g(W+1) - (\b W+W(W-1))g(W)\}| \le \e_0 \|g\| + \e_1 \|\D g\|+\e_2\|\D^2
g\|
$$
for all functions~$g$ satisfying $\mean[W^2|g(W)|]<\infty$,
{where $\Delta g(\cdot):=g(\cdot+1)-g(\cdot)$ and $\Delta^2g(\cdot):=\Delta (\Delta g(\cdot))$.
It} then follows from~\Ref{AX01} that
\beas
d_W(\L(W),\PBD(\a;0,\b,1)) &=& \sup_{f\in\F}|\mean f(W) - \PBD(\a;0,\b,1)(f)|\\
&\le& \e_0 \sup_{f\in\F} \|g_f\| + \e_1\sup_{f\in\F} \|\D g_f\|
  + \e_2\sup_{f\in\F} \|\D^2 g_f\|.
\enas

Hence the major obstacle in applying Stein's method is on getting the right
estimates of the solution $g$ and its first two differences. The necessity of
estimating the second difference of $g$ here comes from the fact that the second order
Palm distribution is used. We summarize the estimates in the following theorem, where $\eta_1\wedge\eta_2:=\min\{\eta_1,\eta_2\}$.

\begin{thm}\label{maintheorem1} The solution $g_f$ to Stein's equation 
\eq\label{eqsteinpbd1}
\a g_f(i+1) -[ \b i+i(i-1)] g_f(i) = f(i) - \PBD(\a;0,\b,1)(f)
 \en
 with convention $g_f(0) = g_f(1)$ satisfies
\bea
\sup_{f \in \F} \| g_f \| &=& \frac{1}{\a} \sum_{k=1}^\infty k \pi_k  \le 
\frac{1}{\b+\frac{2\a}{\a+2\b+2}}\bigwedge\frac{\sqrt{ \a + (1-\b)^2/4 }  + (1-\b)/2 
}{\a}, \label{estimate1}\\
\sup_{f \in \F} \| \D g_f \|&\le&  \leo  \frac{1}{\b+1} +\frac{\a+ 2 }{ \a(\a +2) + 2 \b (\a + \b + 1)  }  \rio \bigwedge \left(  \frac{1}{ \sqrt{\a} } + \frac{1}{ \a} \right) , \label{estimate2}\\
\sup_{f \in \F} \| \D^2 g_f \|&\le&\frac{2}{ \a } +  \frac{1}{ \a + \b ( 1+ \b +\a / 2) }\le\frac{3}{\a}.\label{estimate3}
\ena
\end{thm}

Some care is needed when choosing the parameters in PBD distributions. For example,
instead of using $\PBD(\a;0,\b,1)$, one may prefer to use $\PBD(a;0,1,b)$ with
$\a_i=a$ and $\b_i=i+bi(i-1)$, $i\ge 0$, since{, when $b=0$,} the distribution is
reduced to {the Poisson distribution with mean $a$, denoted by} $\Pn(a)$. In other words, $\PBD(a;0,1,b)$ can be viewed as a
generalization of the Poisson distribution and its Stein equation becomes
\eq\label{eqsteinpbd2}
a\tg_f(i+1)-[i+bi(i-1)]\tg_f(i)=f(i)-\PBD(a;0,1,b)(f).
\en
To relate to \Ref{eqsteinpbd1}, one simply needs to define $a=\a/\b$, $b=1/\b$ and
$\tg_f=\b g_f$. The following proposition is an immediate corollary of
Theorem~\ref{maintheorem1}.

\begin{cor} The solution $\tg_f$ to Stein's equation \Ref{eqsteinpbd2}  satisfies
\bea
\sup_{f \in \F} \| \tg_f \| &\le& \frac{1}{1+\frac{2ab}{a+2b+2}}\bigwedge \frac{\sqrt{
ab + (b-1)^2/4 }  + (b-1)/2  }{ab}, \label{estimate4}\\
\sup_{f \in \F} \| \D \tg_f \|&\le&  \leo  \frac{1}{1+b} + \frac{ a+ 2b }{a^2+2a+2+2b(a+1)}  \rio \bigwedge \left( \frac{1}{ \sqrt{ab}}  + \frac{1}{ a} \right),\label{estimate5}\\
\sup_{f \in \F} \| \D^2 \tg_f \|&\le&\frac{2}{ a } +  \frac{b}{ (a + 1) b+ 1 +a /2 }.\label{estimate6}
\ena
\end{cor}

\begin{Remark}{\rm It is interesting to note that for $\Pn(a)$ approximation, i.e.
$b=0$, Barbour and Xia~(2006) showed that the corresponding estimates of $\tg_f$ for
$\Pn(a)$ are
\bea
\sup_{f\in\F}\|\tg_f\|&=& 1,\label{first}\\
\sup_{f\in\F}\|\D \tg_f\|&\le& 1\wedge\frac{8}{3\sqrt{2ea}}
  \ \le\ 1\wedge\frac{1.1437}{\sqrt{a}},\label{second}\\
\sup_{f\in\F}\|\D^2\tg_f\|&\le&\frac{4}{3}\wedge \frac{2}{ a }.\label{third}
\ena
}\end{Remark}
It is straightforward to check that when $b=0$, estimates \Ref{estimate4} and
\Ref{estimate6} are reduced to \Ref{first} and \Ref{third} respectively, while 
\Ref{estimate5} becomes $\frac{a^2+3a+2}{a^2+2a+2}$, which is not the same as \Ref{second}.
It is an interesting unsolved problem to establish a bound of $\sup_{f \in \F} \|
\D \tg_f \|$ for $\PBD(a;0,1,b)$ which reduces to \Ref{second} when $b=0$.

\ignore{As pointed out in Barbour and Xia~(2006), the differences of $\Pn(a)$ approximation in terms of
the Wasserstein distance typically occur 
at places separated by distances of order~$\sqrt a$, the standard deviation of
$\Pn(a)$. Hence the error estimates of $d_W(\L(W),\Pn(a))$ are expected typically of the order of that for the
total variation distance magnified by a factor of $\sqrt a$. In order to improve on
the accuracy, one needs to match not only the means of $W$ and the approximating
distribution but also their variances. This consideration leads to various modified
Poisson approximations such as shifted (or translated) Poisson and compound Poisson
signed measures where negative `probabilities' are often needed [see Presman~(1983),
Kruopis~(1986), \v Cekanavi\v cius~(1997), Barbour \& Xia~(1999), \v Cekanavi\v cius \& Vaitkus~(2001) and
R\"ollin~(2005)]. PBD distributions seem to be a more natural }

A prototypical example for applying the above estimates is to consider {the Poisson binomial distribution, i.e.} the {distribution of the} sum $W$
of independent Bernoulli random variables $\{X_i : 1 \le i \le n \}$ with
{the} distribution
 $$
\prob(X_i = 1) = 1 - \prob(X_i = 0) = p_i, \ \ \ 1 \le i \le n.
 $$
Denote $ \lambda_r = \sum_{i=1}^n p_i^r${, $\theta_r=\l_r/\l_1$} for $r \ge 1$ and $\l:=\l_1$. 
Barbour and Xia~(2006) proved that if $\l_2$ is an integer, then
\eq
  d_{W}(\L(W),\Pn(\l-\l_2)\ast\delta_{\l_2}) \le  4(\l-\l_2)^{-1}\l_2{=\frac{4\theta_2}{1-\theta_2}}, \label{BX}
\en
 where $*$ denotes convolution and $\delta_{\l_2}$ the point mass at~$\l_2$.
However, when $\l_2$ is not an integer, some minor adjustment is needed to make the
argument work and the error bound will have to be increased by an amount reflecting
the correction. The following theorem shows an impressive improvement on the quality of approximation by a PBD distribution {since the bound in \Ref{application1} below is of order $\theta_3+\theta_2^2$ while the bound in \Ref{BX}  is of order $\theta_2$ only.}

 \begin{thm} With the setup in the preceding paragraph, set
$$
\b = \l^2 \lambda_2^{-1} - 1 - 2 \lambda + 2 \lambda_3 \lambda_2^{-1}, \ \ \ \a
= \b \lambda + \lambda^2 - \lambda_2.
$$
Provided $\l^2 \lambda_2^{-1} - 1 - 2 \lambda\ge0$, we have
\eq\label{application1}
d_W ({\cal L} W, \PBD(\a;0,\b,1) ) \le  3 \theta_3 + \frac{ 6 \theta_2 \l_2 }{  \l  - \lambda_2 - (1 + \theta_2) \theta_2 } .
 \en

\end{thm}

\noindent{\bf Proof.} The expansion of $\ex \B g (W)$ in Brown and Xia~(2001), p.~1390 states that
\bea
\mean [\B g(W)]&=&-\b\sum_{i=1}^np_i^3\mean[\D^2g(W^i+1)]+\sum_{j\ne i}p_i^2p_j^2(1-p_i-p_j)\mean[\D^2g(W^{ij}+2)]\nonumber\\
&&+\sum_{j\ne i}p_ip_j(p_i+p_j)(1-p_i)(1-p_j)\mean[\D^2g(W^{ij}+1)],\label{BrownXiaexapansion}
\ena
where $W^i=W-X_i$ and $W^{ij}=W-X_i-X_j$. 
{Therefore,
\beas
&&d_W(\L(W),\PBD(\a;0,\b,1)) = \sup_{f\in\F}|\mean f(W) - \PBD(\a;0,\b,1)(f)|\\
&&=\sup_{f\in\F} |\mean [\B g_f(W)]|\\
&&\le \sup_{f \in \F} \| \D^2 g_f \|\left\{\b\sum_{i=1}^np_i^3+\sum_{j\ne i}p_i^2p_j^2(1-p_i-p_j)+\sum_{j\ne i}p_ip_j(p_i+p_j)(1-p_i)(1-p_j)\right\}\\
&&\le (\b \lambda_3 + 2 \l \l_2) \sup_{f \in \F} \| \D^2 g_f \|\\
&&\le3 (\b \lambda_3 + 2 \l \l_2)/\a ,
\enas}
{where we applied \Ref{eqsteinpbd1} in the second equation, \Ref{BrownXiaexapansion} in the first inequality and \Ref{estimate3} in the last inequality.}
Since $\b \ge \l^2 \lambda_2^{-1} - 1 - 2 \lambda$ and $\l^2\ge \l_2$, we have $\a {= \b \lambda + \lambda^2 - \lambda_2}\ge  ( \l^2 \lambda_2^{-1} - 1 - 2 \lambda ) \l + \l^2 - \l_2 = \l^3 \l_2^{-1} - \l - \l^2 - \l_2$ and $\a\ge \b\l$. Therefore, the above upper bound can be further estimated as follows:
 $$
\frac{3\b \l_3}{\a} + \frac{6 \l \l_2}{\a} \le {\frac{3\b\l_3}{\b\l}}+ \frac{6 \l \l_2}{\l^3  \l_2^{-1} - \l - \l^2 - \l_2} = 3 \theta_3 + \frac{ 6 \theta_2 \l_2 }{  \l  - \lambda_2 - (1 + \theta_2) \theta_2 } .\  \qed
 $$
  
\begin{Remark}{\rm If one uses \Ref{estimate2} and argues as in the proof of Theorem~3.1 of Brown and Xia~(2001), the following estimate can be established: }
$$ 
d_W ({\cal L} W, \PBD(\a;0,\b,1) ) \le \leo  \frac{\b \lambda_3}{\sigma_1} + \frac{2
\lambda \lambda_2}{ \sigma_2} \rio \leo  \frac{1}{\sqrt{\a}} + \frac{1}{\a} \rio \asymp O\left(\sqrt{\b}\theta_3{+\sqrt{\l}\theta_2^{1.5}}\right),
 $$
{\rm where $\sigma_k=\sqrt{\sum_{i=k+1}^n \rho_i}$, $\rho_i$ is the $i$th largest value of $p_1(1-p_1)$, $p_2(1-p_2)$, $\cdots$, $p_n(1-p_n)$. {The bound here is not as good as the bound in \Ref{application1} since the latter is of order $\theta_3+\theta_2^2$ while $\sqrt{\l}$ and $\sqrt{\b}\asymp \l/\sqrt{\l_2}$ are usually large}.}
\end{Remark}

\begin{Remark}{\rm Both $\PBD(\a;0,\b,1)$ and shifted Poisson use two parameters,
however, the bound in \Ref{application1} is significantly better than \Ref{BX}. In
most applications, we can not expect $\l_2$ to be an integer which means some
correction in shifted Poisson approximation is unavoidable, while the PBD
approximation does not have this inconvenience.}
\end{Remark}

\section{The proofs.}
\setcounter{equation}{0}

We follow the idea of Brown \& Xia~(2001) by setting $g_f(i)=h_f(i)-h_f(i-1)$ for $i\ge 1$ so
that
\eq \label{generator}
\B g_f(i) = \a_i [h_f(i+1)-h_f(i)] + \b_i [h_f(i-1)-h_f(i)] :=\A h_f(i),\ i\ge 0,
 \en
where $\A$ is the generator of birth-death Markov chains with birth rates $\{\a_i:\
i\ge0\}$ and death rates $\{\b_i:\ i\ge1\}$. In this setup, Stein's equation
\Ref{Steinequ} becomes
$$\A h_f(i)=f(i)-\pi(f) , $$
and its solution can be explicitly expressed as
\eq
h_f(i)=-\int_0^\infty\mean [f(Z_i(t))-\pi(f)]dt,\label{solution01}
\en 
where $\{Z_i(t):\ t\ge0\}$ is the birth-death process with generator $\A$ and initial state $i\ge 0$.

The advantage in estimating Stein's constants in Poisson approximation is that the couplings are relatively trivial and many computations are achievable [see Barbour and Xia~(2006)]. For the $\PBD(\a;0,\b,1)$ setting, we follow the main ideas in Barbour and Xia~(2006) but use immigration-death particle systems in $\{ 0,1 \}^\Z$ to realize the couplings. Due to the generality of the model, our estimations are naturally more technical.

\ignore{The main ideas of the proofs are similar to those in Barbour and Xia~(2006).
However, given the generality of the model, there are a lot of striking differences
between $\PBD(\a;0,\b,1)$ and a Poisson distribution since in the
latter case, the couplings are relatively trivial and many computations are
achievable. }

Let $g_k$ be the solution of (\ref{Steinequ}) for $f=\bone_k$, then Lemma~2.3 of
Brown \& Xia~(2001) states that
 \eq \label{gk}
g_k (i) = \frac{ \pi_k \bar{F} (i) }{ \a_{i-1} \pi_{i-1} }, \  {\rm if} \ k< i ; \ \ \
 g_k (i) = \frac{ - \pi_k F(i-1) } { \b_i \pi_i },  \  {\rm if} \ k \ge  i,
 \en
where $\pi$ is defined in \Ref{equilibriumdis}, $F(i) = \sum_{j \le i} \pi_j$ and
$\bar{F}(i)  = \sum_{j \ge i} \pi_j$.

\subsection{The proof of \Ref{estimate1}.} \label{PBDg}

Noting that $g_f(0):=g_f(1)$ and  $\sup_{f \in \F}  \| g_f\|=\sup_{i\ge
1}\sup_{f\in\F}|g_f(i)|$, we will show in Lemma~\ref{Lemma11} that
$$\sup_{f\in\F}|g_f(i)|=g_{f_1} (i), $$
where
$$f_{1}(k)=-k,\ k\in\Z_+.
$$
We then show in Lemma~\ref{lem2} that $g_{f_{1}}(i)$ is a decreasing function in
$i$ and reduce the estimate to 
$$
\sup_{f \in \F} \| g_f \| = g_{f_{1}}(1) = \frac{1}{\a} \sum_{k=1}^\infty k \pi_k ,
 $$
where the last equation comes from solving the Stein equation \Ref{Steinequ} with $i=0$.  The
estimate $(\b+\frac{2\a}{\a+2\b+2})^{-1}$ comes from \Ref{epiestimate} below. Now we concentrate on the other part of the upper bound in \Ref{estimate1}.

Noting that {$\a_k=\a$, $\b_k=\b k+k(k-1)$ and the balance equation $\a_k\pi_k=\b_{k+1}\pi_{k+1}$ holds for all $k\in\Z_+$, we have}
 \beas
\sum_{k=1}^\infty k \pi_k
 & \le &
\left( \sum_{k=1}^\infty k^2 \pi_k  \right)^{1/2} {=} \left( \sum_{k=1}^\infty
k(k-1) \pi_k + \sum_{k=1}^\infty k \pi_k  \right)^{1/2}
 \\ & = &
\left( \sum_{k=1}^\infty (\b_k - \b k) \pi_k + \sum_{k=1}^\infty k \pi_k 
\right)^{1/2}
=  \left( \sum_{k={1}}^\infty \a_{{k-1}} \pi_{{k-1}}+ (1-\b) \sum_{k=1}^\infty k \pi_k 
\right)^{1/2}
  \\ & = &
\left( \a + (1-\b) \sum_{k=1}^\infty k \pi_k  \right)^{1/2} ,
     \end{eqnarray*}
 {hence}
 $$
\sum_{k=1}^\infty k \pi_k  \le \big( \a + (1-\b)^2/4 \big)^{1/2} + (1-\b)/2.
 $$
This implies
 $$
\sup_{f \in \F} \| g_f\|  \le \frac{\sqrt{ \a + (1-\b)^2/4 }  + (1-\b)/2  }{\a} ,
 $$
completing the proof of \Ref{estimate1}.

\begin{lem} \label{Lemma11} For each $i\ge 1$, we have
$$\sup_{f\in\F}|g_f(i)|=g_{f_{1}}(i).
$$
\end{lem}

\noindent{\bf Proof.} Since the solution of (\ref{Steinequ}) does not change when
$f$ is shifted by a constant, it suffices to consider $f\in\F$ with $f(i)=0$, that
is,
 $$
\sup_{f\in\F}|g_f(i)|= \sup_{f \in \F : f (i) = 0} \left| \sum_{k=0}^\infty f(k) g_k
(i) \right|.
 $$
From (\ref{gk}), we know that $g_k (i)$ is positive for $k<i$ and negative for $k
\ge i$. Given that $f \in \F$ with $f(i) = 0$, it is easy to see that $f(k) \le
|k-i|$ and 
 $$
|g_f(i)| \le \sum_{k : k <i}^\infty |i-k| g_k (i) + \sum_{k : k > i}^\infty (-
|i-k|) g_k (i),
 $$
with the equality attained at $f (k) = i-k$, or equivalently, at $f = f_1$.  \qed
     
The following lemma summarizes the major properties of $g_{f_{1}}$.
 
\begin{lem}\label{lem2} \label{delta2ti}
The functions $g_{f_{1}}$ and $\D^2 g_{f_1}$ are positive and decreasing in $i\ge 1$, and $\D g_{f_1}$ is increasing in $i \ge 1$.
\end{lem}

\noindent{\bf Proof.} We will exploit a probabilistic interpretation of $g_{f_1}$ by constructing a coupling of $\{ Z_i 
(t) : t  \ge 0 \}$ for all $i \in \Z_+$ as follows [cf Barbour~(1988) and Barbour \&
Brown~(1992)]. Consider a particle system on the site space $\Z$, the set of all
integers. For each configuration $\xi \in \{ 0,1 \}^\Z$, set $\xi_n = 1$ if position
$n$ is occupied by a particle, and 0 if position $n$ is vacant.  Let $\{\xi (t) : t
\ge 0 \}$ be a Markov process with state space $\{ 0,1 \}^\Z$ and initial state $\xi_n
(0) = \bone_{ \{ n \ge 1 \} }$, which evolves as follows:
\begin{enumerate}
\item particles immigrate to the system with rate $\a$ and a new immigrant takes the
closest vacant site to the left of site 1;
\item each existing particle suicides with rate $\b$;
\item the particle at site $i$ kills the particle at site $j (>i)$ with rate 2;
\item all above evolution rules apply independently.
\end{enumerate}
Fix $i \ge 0$. The evolution of $\sum_{n = - \infty}^i \xi_n (t)$ does not depend on
the particles at the sites to the right of $i$, i.e. $\{ \xi_n (t) : t \ge 0, n > i
\}$. It is straightforward to check that $\{ \sum_{n = - \infty}^i \xi_n (t) : t \ge 0
\}$ is a birth-death process with birth rate $\a_k = \a$ and death rate
$\b_k = \b k + k(k-1)$, which means it has the same distribution as $\{ Z_i (t): t \ge 0 \}:=\bZ_i$. 
Take \begin{equation}\label{expressionZ} Z_i (t) = \sum_{n= -\infty}^i \xi_n (t) \end{equation} from now on.
Let 
\eq\label{coupletime}
T_i : = \inf \{ t \ge 0 : \xi_i (t) = 0 \}
\en
 be the first time that position $i$ becomes vacant.

Noting that $g_f(i)=h_f(i)-h_f(i-1)$ with $h_f$ in
\Ref{solution01}, we have that, for any $i \ge 1$,
 \bea
 g_{f_{1}}(i) 
 & = & 
h_{f_{1}}(i) - h_{f_{1}}(i-1) =\mean \int_0^\infty ( Z_i (t) -
Z_{i-1} (t) ) d t \nonumber
 \\ &  = & 
\mean\int_0^\infty \xi_i (t) d t = \mean
T_i=\int_0^\infty\prob(T_i>t)dt,\label{Xiaaddagain2}
 \ena
 which implies that $g_{f_1}(i)\ge 0$.
Moreover, for any $i \ge 1$, 
 \beas
\prob(T_i > t)
 = 
\mean\left[ \prob \big( T_i > t | \bZ_{i-1} \big)\right]
  = \mean \exp \left\{ - \int_0^t \left(\b + 2 Z_{i-1}(s)\right) d s \right\},
 \enas
which is decreasing in $i$ for any $t$ since $Z_{i-1}(s) \le Z_i (s)$ for any $i \ge 1$ and $s \ge 0$. Hence $g_{f_1} (i)$ decreases in $i$.

Integration by parts ensures
 \begin{eqnarray*}
\mean (e^{- 2 (t \wedge T_i ) }  |\bZ_{i-1} )
 & =&
 \int_0^\infty e^{- 2 (s\wedge t)} d \big( - \prob ( T_i > s |
\bZ_{i-1} ) \big)
  \\ & =&
1 -  \int_0^t 2 e^{- 2 s}  \prob (T_i  > s |\bZ_{i-1} ) d s.
  \end{eqnarray*}
Consequently, the increasing property of $\D g_{f_1}$ follows from
$$\D g_{f_1}(i)=\int_0^\infty (\prob(T_{i+1}>t)-\prob(T_i>t))dt$$
and
 \beas
 \prob (T_i > t) - \prob(T_{i+1} > t)
  & = &
\mean \leo  e^{- \int_0^t (\b + 2 Z_{i-1} (s)) d s} \leo  1 - \mean ( e^{- 2 (t \wedge T_i )
}  |\bZ_{i-1} ) \rio  \rio \nonumber
 \\ & = &
\mean \leo  e^{- \int_0^t (\b + 2 Z_{i-1} (s)) d s}   \int_0^t 2 e^{- 2 s}  \prob (T_i >
s |\bZ_{i-1}  ) d s \rio \nonumber
  \\ & = & 
\mean \leo  e^{- \int_0^t (\b + 2 Z_{i-1} (s)) d s}  \int_0^t 2 e^{- 2 s} e^{-\int_0^s (\b+ 2 Z_{i-1} (u)) d u} d s \rio ,
 \enas
 which is decreasing in $i$.
Using the above equation again,  
we conclude that
 \begin{eqnarray*}
 &&
\prob (T_i > t) - 2 \prob(T_{i+1} > t) + \prob(T_{i+2} > t)
  \\  & = &
\mean \leo  e^{- \int_0^t (\b + 2 Z_{i-1} (s)) d s}   \int_0^t 2 e^{- 2 s}  \left[
e^{-\int_0^s (\b + 2 Z_{i-1}(u)) d u } 
- \mean \big(  e^{- 2 (t \wedge T_i )} e^{-\int_0^s( \b + 2 Z_i (u)) d u } |\bZ_{i-1}
\big) \right] d s  \rio 
 \\ & = &
\mean \leo  e^{- \int_0^t( \b + 2 Z_{i-1} (s)) d s}   \int_0^t 2 e^{- 2 s} e^{-\int_0^s
(\b + 2 Z_{i-1}(u)) d u }  \left[ 1
- \mean \big(  e^{- 2 (t \wedge T_i ) - 2 (s \wedge T_i ) } |\bZ_{i-1} \big) \right]
d s  \rio .
   \end{eqnarray*}
The decreasing property of $\D^2 g_{f_1}$ follows from the facts that 
$$\D^2g_{f_1}(i)=\int_0^\infty(\prob(T_{i+2}>t)-2\prob(T_{i+1}>t)+\prob(T_i>t))dt,$$
$Z_{i-1}(s)$ is increasing in $i$ and
\beas
1-\mean \left(e^{-2(t\wedge T_i)-2(s\wedge T_i)}|\bZ_{i-1}\right)&=&\int_0^\infty \prob(T_i>u|\bZ_{i-1})d\left(-e^{-2(t\wedge u)-2(s\wedge u)}\right)\\
&=&\int_0^\infty e^{-\int_0^u(\b+2Z_{i-1}(v))dv}d\left(-e^{-2(t\wedge u)-2(s\wedge u)}\right),
\enas
which is decreasing in $i$. Finally, $\D^2g_{f_1}$ is positive because $\D g_{f_1}$ is increasing. \qed

\subsection{The proof of \Ref{estimate2}.} \label{PBDdeltag}

Arguing as in subsection~\ref{PBDg}, we conclude that
 $$
\sup_{f \in \F}   \|\D g_f \| =\sup_{i\in\Z_+} \sup_{f \in \F}  \left|
\sum_{k=0}^\infty f(k) \D g_k (i) \right|=\sup_{i\in\Z_+} \sup_{f \in \F : f(i) = 0}
 \left| \sum_{k=0}^\infty f(k) \D g_k (i) \right|.
 $$
Since $\delta_i = (\b_{i+1} - \b_i ) - (\a_{i+1} - \a_i) = \b + 2 i$ is positive, 
Lemma~2.4 {with condition (C4)} in Brown \& Xia~(2001) states that $ \D g_k (i) $ is negative whenever $k
\neq i$, which ensures that, for each $f\in\F$ with $f(i)=0$,
 \bea
&&|\D g_f(i) |\le
\sum_{k=0}^\infty |i-k| (- \D g_k (i))   \label{Si1}
 \\ & =&
\sum_{k: k<i} (i-k) \pi_k \leo  \frac{\bar{F} (i)}{\a \pi_{i-1}} - \frac{ \bar{F} (i+1)
}{\a \pi_i } \rio + \sum_{k: k>i} (k-i) \pi_k \left( \frac{ F (i) }{\b_{i+1} \pi_{i+1}
} - \frac{ F (i-1) }{\b_i \pi_i } \right),  \nonumber    \ena
with the equality in \Ref{Si1} reached when
$$f(k)=-|i-k|={:}f_{i2}(k).$$
However, since $\sum_{k\in\Z_+}\D g_k(i)=0$, we have
  \bea
\sum_{k: k \neq i} \big( - \D g_k (i) \big) \nonumber
&=&
\sum_{k: k<i}  \pi_k \leo  \frac{\bar{F} (i)}{\a \pi_{i-1}} - \frac{ \bar{F} (i+1) }{\a
\pi_i } \rio + \sum_{k: k>i}  \pi_k \left( \frac{ F (i) }{\b_{i+1} \pi_{i+1} } -
\frac{ F (i-1) }{\b_i \pi_i } \right)\\
&=& \D g_i (i) = \frac{ \bar{F}(i+1) }{ \a_i} + \frac{ F(i-1) }{ \b_i } \le
\frac{1}{\a_i} \wedge \frac{1}{\b_i},\label{Xiaaddagain1}
  \ena
where the last inequality is due to Theorem 2.10 in Brown \& Xia~(2001). Now we
compare the coefficients $ \sum_{k: k<i} (i-k) \pi_i $ with $\sum_{k : k < i} \pi_k
= F(i-1)$, as well as $\sum_{k: k > i} (k-i) \pi_k$ with $\bar{F} (i+1)$. For the
former, we use $i-k \le i$. Apropos of the latter, we apply the Cauchy-Bunyakovskii
inequality to get
 \beas
\sum_{k: k>i} (k-i) \pi_k
 & = &
\bar{F} (i+1) + \sum_{k: k>i+1} (k-i-1) \pi_k
 \\ & \le &
\bar{F} (i+1) + \left( \sum_{k: k>i+1} (k-i-1)^2 \pi_k \right)^{1/2} \left( \sum_{k:
k>i+1} \pi_k \right)^{1/2}
 \\ & \le &
\bar{F} (i+1) + \left( \sum_{k: k>i+1} \b_k \pi_k \right)^{1/2} \sqrt{\bar{F} (i+1)}
 \\& = &
\bar{F} (i+1) + \left( \sum_{k: k>i} \a_k \pi_k \right)^{1/2} \sqrt{ \bar{F} (i+1) }
 \\ & = &
(1 + \sqrt{\a} )  \bar{F} (i+1).
\end{eqnarray*}
It now follows from \Ref{Xiaaddagain1} that
\beas
|\D g_{f_{i2}}(i)|
& \le &
i F(i-1) \left( \frac{ \bar{F} (i) }{\a \pi_{i-1} } - \frac{ \bar{F} (i+1) }{\a
\pi_i } \right) + ( \sqrt{\a} + 1) \bar{F} (i+1)  \left( \frac{ F (i) }{\b_{i+1}
\pi_{i+1} } - \frac{ F (i-1) }{\b_i \pi_i } \right)
 \\ & \le &
 \max \{ i,  \sqrt{\a} + 1 \} \sum_{k : k \neq i} (- \D g_k (i))
 \le
\max \{ i,  \sqrt{\a} + 1 \}  \left( \frac{1}{\a_i} \wedge \frac{1}{\b_i} \right).
 \end{eqnarray*}
Hence, for $i \le \sqrt\a + 1$,
  $$
|\D g_{f_{i2}}(i)|  \le \frac{\sqrt{\a} + 1}{\a},
  $$
and for $i > \sqrt\a + 1$,
  $$
|\D g_{f_{i2}}(i)| \le \frac{i}{\b_i} \le \frac{1}{\b + i - 1} < \frac{1}{\b +
\sqrt{\a}} < \frac{1}{\sqrt{\a}},
  $$
  which in turn ensure
 $$
\sup_{f \in \F} \| \D g_f \| \le 1 / \sqrt \a + 1 / \a.
 $$

On the other hand, by \Ref{solution01}, \Ref{expressionZ} and \Ref{coupletime},
\begin{eqnarray*}
 \D g_{f_{i2}} (i) 
 & = & 
\ex  \int_0^\infty (| Z_{i+1} (t) - i | - 2 | Z_i (t) - i | + | Z_{i-1} (t) - i | )d t 
  \\  & = &
\ex  \int_0^\infty (|  Z_{i-1} (t) - i + 1_{T_i > t} + 1_{T_{i+1} > t} | - 2 | Z_{i-1} (t) - i + 1_{T_i > t} | + | Z_{i-1} (t) - i |) d t  
  \\  & = &
\ex  \int_0^\infty \leo  |  Z_{i-1} (t) - i + 2 | - 2 | Z_{i-1} (t) - i + 1| + | Z_{i-1} (t) - i | \rio 1_{T_i, T_{i+1} > t} d t 
  \\ &&
+ \ex  \int_0^\infty \leo  | Z_{i-1} (t) - i| -  |  Z_{i-1} (t) - i + 1 |\rio \leo  1_{T_i > t, T_{i+1} \le t} - 1_{T_i \le t, T_{i+1} > t} \rio d t 
  \\  & = &
2 \int_0^\infty \prob \leo  Z_{i-1}  (t)  = i - 1, T_i , T_{i+1} > t  \rio  d t ,
 \\ && 
+  \ex \int_0^\infty \leo  1_{Z_{i-1}(t) \le i-1} - 1_{Z_{i-1} (t)\ge i} \rio  \leo   1_{T_i > t} - 1_{T_{i+1} > t} \rio d t  .
  \end{eqnarray*}
Now
  \begin{eqnarray*}
   &&
2\int_0^\infty \prob \leo  Z_{i-1}  (t)  = i - 1, T_i, T_{i+1} > t  \rio  d t
\le 2\int_0^\infty  \prob \leo T_i, T_{i+1} > t  \rio  d t
 \\ & = &
2\ex \int_0^\infty 1_{  T_i > t } e^{- \int_0^t (\beta + 2 Z_{i-1} (s) + 2)ds}  d t=2
\ex \int_0^\infty e^{- 2 \int_0^t (\beta + 2 Z_{i-1} (s) + 1)ds}  d t
 \\ & \le &
2\int_0^\infty e^{- 2 ( \beta + 1)t}  d t = \frac{1}{\beta+1} ,
 \end{eqnarray*} 
and, by \Ref{Xiaaddagain2} and Lemma~\ref{delta2ti},
 \begin{eqnarray*}
 &&
\left| \ex \int_0^\infty \leo  1_{Z_{i-1}(t) \le i-1} - 1_{Z_{i-1} (t)\ge i} \rio  \leo   1_{T_i > t} - 1_{T_{i+1} > t} \rio d t  \right| 
 \\ & = &
\left| \ex  \int_0^\infty \leo  1_{Z_{i-1}(t) \le i-1} - 1_{Z_{i-1} (t)\ge i} \rio  \leo   \prob ( T_i > t | \bZ_{i-1} ) - \prob ( T_{i+1} > t |\bZ_{i-1} ) \rio d t \right|
  \\ & \le &
\ex  \int_0^\infty \left| 1_{Z_{i-1}(t) \le i-1} - 1_{Z_{i-1} (t)\ge i}   \right|  \big( \prob ( T_i > t |\bZ_{i-1} ) - \prob ( T_{i+1} > t |\bZ_{i-1} ) \big) d t
  \\ & \le &
\ex  \int_0^\infty ( \prob ( T_i > t |\bZ_{i-1} ) - \prob ( T_{i+1} > t |\bZ_{i-1} ) ) d t
 \\ & = &
\ex  T_i - \ex  T_{i+1} =-\D g_{f_1}(i)\le - \D g_{f_1} (1) .
 \end{eqnarray*} 
 Solving the Stein equation \Ref{Steinequ} with $i=0$ and 1 gives
 \eq
 g_{f_1}(1)=\frac{\mean \pi}{\a},\ 
g_{f_{1}}(2) = \frac{1}{\a} \leo   \b_1 g_{f_1}(1) + \mean \pi - 1 \rio  = \frac{ ( \a + \b ) \mean \pi
- \a  }{\a^2} ,\label{Xiaaddagain3}
 \en
so $ - \D g_{f_1} (1) =\frac{1}{\a}\left(1-\frac{\b}{\a}\mean \pi\right)$.
Collecting these estimates yields
 \begin{eqnarray*}
\D g_{f_{i2}} (i) 
 & \le &
\frac{1}{\b+1} + \frac{1}{\a}  \leo  1 - \frac{\b \ex \pi}{\a} \rio 
 \le 
\frac{1}{\b+1} + \frac{\a+ 2 }{ \a(\a +2) + 2 \b (\a + \b + 1)  }   , 
 \end{eqnarray*} 
where the last inequality is from \Ref{epiestimate} below. \qed

\subsection{The proof of \Ref{estimate3}.} \label{PBDdeltatwog}

As in Brown \& Xia~(2001), we define $\tau_i^+=\inf\{t:\ Z_i(t)=i+1\}$ for $i\ge 0$
and $\tau_i^-=\inf\{t:\ Z_i(t)=i-1\}$ for $i\ge1$. For convenience, we
set $\tau_0^-=\infty$, $e_i^+=\mean \tau_i^+$ and $e_i^-=\mean \tau_i^-$ for all
$i\ge0$. It has been proved in Brown \& Xia~(2001, p.~1378) that
$$e_i^+ = \frac{F(i)}{\a \pi_i}\mbox{ for }i \ge 0,\ e_i^- = \frac{\bar{F} (i)}{\a
\pi_{i-1}}\mbox{ for }i\ge 1$$
and
\eq
g_j(i)=\left\{
\begin{array}{ll}
-\pi_je_{i-1}^+,&i\le j,\\
\pi_je_i^-,&i\ge j+1.
\end{array}
\right.
\label{Xiaaddagain4}
\en
We begin the proof with two technical lemmas.

\begin{lem}\label{convex} For $i \ge 0$,
$e_{i+2}^+-2e_{i+1}^++e_i^+\ge0$ and for $i\ge 1$, $e_{i+2}^--2e_{i+1}^-+e_i^-\ge 0 $.
 \end{lem}

\noindent{\bf Proof.} For $i \ge 0$,  
 \begin{eqnarray*}
  &&
\frac{F(i+2)}{\a \pi_{i+2}} - 2 \frac{F(i+1)}{\a \pi_{i+1}} + \frac{F(i)}{\a \pi_i}
 \\ & = &
\frac{1}{\a^3 \pi_i} \leo  \b_{i+2} \b_{i+1} F(i+2) - 2 \b_{i+1} \a F(i+1) + \a^2 F(i)
\rio 
 \\ & =&
\frac{1}{\a^3 \pi_i} \leo  \b_{i+2} \b_{i+1} \sum_{k=0}^{i+2} \pi_k - 2 \b_{i+1}
\sum_{k=1}^{i+2} \b_k \pi_k + \sum_{k=2}^{i+2} \b_k \b_{k-1} \pi_k \rio 
\\ & =&
\frac{1}{\a^3 \pi_i} \leo   \b_{i+2} \b_{i+1} \pi_0 + (\b_{i+2} \b_{i+1}  - 2 \b_{i+1} 
\b_1 ) \pi_1 + \sum_{k=2}^{i+1} (\b_{i+2} \b_{i+1} - 2  \b_{i+1}  \b_k +  \b_k
\b_{k-1} ) \pi_k  \rio .
 \end{eqnarray*}
Direct computation shows that $\b_{i+2}\ge 2\b_1$ and 
 \beas
&&\b_{i+2} \b_{i+1} - 2  \b_{i+1}  \b_k +  \b_k \b_{k-1}\\
&=&\b_{i+1}^2+(\b+2(i+1))\b_{i+1}-2\b_{i+1}\b_k+\b_k^2-\b_k(\b+2(k-1))\\
&=&(\b_{i+1}-\b_k)^2+\b_{i+1}(\b+2(i+1))-\b_k(\b+2(k-1)) \ge 0
\enas
for $2 \le k \le i+1$, which implies the claim for $\{e_i^+:\ i\ge 0\}$. 
 
For $i \ge 1$,
 \begin{eqnarray*}
  &&
\frac{\bar{F}(i+2)}{\a \pi_{i+1}} - 2 \frac{\bar{F}(i+1)}{\a \pi_i} +
\frac{\bar{F}(i)}{\a \pi_{i-1}}
 \\ & = &
\frac{1}{\a^3 \pi_{i-1}} \leo  \b_{i+1} \b_i \bar{F} (i+2) - 2 \b_i  \a \bar{F} (i+1) +
\a^2 \bar{F} (i) \rio  
   \\ & =&
\frac{1}{\a^3 \pi_{i-1}} \sum_{k=i+2}^\infty \leo  \b_{i+1} \b_i - 2 \b_i \b_k +  \b_k
\b_{k-1}  \rio \pi_k,
 \end{eqnarray*}
 so it remains to show $\b_{i+1} \b_i - 2 \b_i \b_k +  \b_k
\b_{k-1}\ge 0$ for $k\ge i+2$.
Since 
$$\frac{\b_{k-1}}{\b_k}=\frac{(k-1)(\b+k-2)}{k(\b+k-1)}=\left(1-\frac{1}{k}\right)\left(1-\frac{1}{\b+k-1}\right),$$
which is increasing in $k$, we can conclude that for $k\ge i+2$, $\b_{k-1}\b_{i+1}\ge \b_k\b_i$, which in turn gives 
 $$
\b_{i+1} \b_i - 2 \b_i \b_k +  \b_k \b_{k-1}\ge 2\sqrt{\b_{i+1}\b_i\b_k\b_{k-1}}-2\b_i\b_k \ge 0,
 $$
completing the proof of the claim for $\{e_i^-:\ i\ge 1\}$. \qed

  \begin{lem}\label{lemdelta2g1} We have
 $$
\D^2 g_{f_{1}}(1) 
 \le
\frac{2}{ \a \b + 2 ( \a + \b + \b^2) } .
$$
 \end{lem}

\noindent{\bf Proof.} Using \Ref{Xiaaddagain3} and the Stein equation \Ref{Steinequ} with $i=2$, we
have
\begin{eqnarray*}
g_{f_{1}}(3)
 & = & 
 \frac{1}{\a} (\b_2 g_{f_1}(2) + \mean \pi - 2)   \\ & = & 
 \frac{1}{\a^3} \leo    \leo   2 (\b + 1 ) (\a + \b)  + \a^2 \rio \mean \pi - 2 \a  (\a + \b
+ 1 )  \rio .
  \end{eqnarray*}
 Using Lemma~\ref{lem2}, we have $0 \le g_{f_{1}}(3) \le
g_{f_{1}}(2)$, which ensures
\eq \label{epiestimate}\frac{2 \a  (\a + \b + 1 ) }{ 2 (\b + 1 ) (\a + \b)  + \a^2} 
\le  \mean \pi \le \frac{\a(\a+2\b+2)}{\b(\a+2\b+2) + 2\a }.
 \en
Therefore,
 \begin{eqnarray*}
\D^2 g_{f_1}(1) 
  = 
 \frac{ 2 (\a + \b + \b^2 ) \mean \pi  - 2 \a ( 1 + \b )  }{\a^3} 
   \le
\frac{2}{ \a \b + 2 ( \a + \b + \b^2) }.
 \end{eqnarray*}
This completes the proof of the lemma. \qed
 
\noindent{\bf The proof of \Ref{estimate3}.} It follows from the definition of 
$g_f(0)=g_f(1)$, \Ref{Si1} and Lemma \ref{lem2}  that
\beas
\left|\D^2 g_f(0)\right|&=&\left|\D g_f(1)\right|\le \D g_{f_{12}}(1)=\D
g_{f_{1}}(1)-2\D g_0(1)
 \le -2\D g_0(1)\le\frac{2}{\a},
\enas
where the last inequality is due to Theorem~2.10 of Brown \& Xia~(2001), so it remains to estimate $\D^2 g_f(i)$ for $i\ge 1$.
In fact, we have from \Ref{Xiaaddagain4} that
 $$\D^2 g_k (i)=\left\{
\begin{array}{ll}
-\pi_k\left(e_{i+1}^+-2e_i^++e_{i-1}^+\right),&\mbox{ for }k\ge i+2,\\
\pi_{i+1}\left(e_{i+2}^-+2e_i^+-e_{i-1}^+\right),&\mbox{ for }k=i+1,\\
\pi_k\left(e_{i+2}^--2e_{i+1}^-+e_i^-\right),&\mbox{ for }k\le i-1.
\end{array}\right.$$
Lemma~2.4 of Brown \& Xia~(2001) ensures that $e_i^+$ is increasing in $i$, so $e_{i+2}^-+2e_i^+-e_{i-1}^+\ge 0$, which, together with Lemma~\ref{convex}, implies
\eq
\D^2 g_k (i)\left\{ \begin{array}{ll} \ge 0, &\mbox{ if } k \le i-1\mbox{ or }k=i+1, \\ \le 0,
&\mbox{ if } k \ge i+2. \end{array} \right.\label{Xiaaddagain5}
\en
Replacing $f$ by $-f$ if necessary, it suffices to give an upper bound for $\D^2g_f(i)$. Again, we fix $i\ge 0$. Using $f-f(i)$ instead of $f$ if necessary, we may assume
$f(i)=0$, then
\beas
\D^2 g_f (i)&=&\sum_{k \neq i} f (k) \D^2 g_k (i)\\
&=&\sum_{k < i} f (k) \D^2 g_k (i)+\sum_{k>i+1} [f (k)-f(i+1)] \D^2 g_k (i)
+f(i+1)\sum_{k\ge i+1}\D^2g_k(i)\\
&\le&\sum_{k < i} (i-k) \D^2 g_k (i)+\sum_{k>i+1} [(i+1)-k] \D^2 g_k (i)
+f(i+1)\sum_{k\ge i+1}\D^2g_k(i),
 \enas
 with the equality reached when 
 $$f(k)=(i-k)\bone_{\{k\le i\}}+[(i+1)-k+f(i+1)]\bone_{\{k\ge i+1\}}:=f_{i3}(k).$$
If $\sum_{k\ge i+1}\D^2g_k(i)<0$, then
$$\D^2 g_{f_{i3}} (i)\le \sum_{k < i} (i-k) \D^2 g_k (i)+\sum_{k>i+1} [(i+1)-k] \D^2
g_k (i)
-\sum_{k\ge i+1}\D^2g_k(i),$$
with the equality attained when $f_{i3}(k)=(i-k)\bone_{\{k\le i\}}+(i-k)\bone_{\{k\ge
i+1\}}=i-k=f_{1}(k) + i$, which is a shifted function of $f_1$. On the other hand,
if $\sum_{k\ge i+1}\D^2g_k(i)\ge0$, then
$$\D^2 g_{f_{i3}} (i)\le \sum_{k < i} (i-k) \D^2 g_k (i)+\sum_{k>i+1} [(i+1)-k] \D^2
g_k (i)
+\sum_{k\ge i+1}\D^2g_k(i) , $$
with the equality achieved at $f_{i3}(k)=(i-k)\bone_{\{k\le i\}}+(i+2-k)\bone_{\{k\ge
i+1\}}:=f_{i4}(k) . $

By Lemma~\ref{lem2},  
$$
\D^2 g_{f_{1}} (i) \le  \D^2 g_{f_{1}} (1).
 $$ 
On the other hand, since $f_{i4}=2\bone_{\{k\ge i+1\}}+f_{1}+i$, we get from \Ref{Xiaaddagain5} and 
\Ref{gk} that
 \begin{eqnarray*}
\D^2 g_{f_{i4}} (i) 
 & = & 
\D^2 g_{f_{1}} (i) + 2 \sum_{k=i+1}^\infty \D^2 g_{k} (i)
 \\ & \le & 
\D^2 g_{f_{1}} (1) + 2 \D^2 g_{i+1} (i) =  \D^2 g_{f_{1}} (1) +  2 \big( g_{i+1}
(i+2) - 2 g_{i+1} (i+1) + g_{i+1} (i) \big)
 \\ & = & 
\D^2 g_{f_{1}} (1) +  2 \pi_{i+1} \leo  \frac{\bar{F} (i+2)}{\a \pi_{i+1}} + 2
\frac{F(i)}{\a \pi_i} - \frac{F(i-1)}{\a \pi_{i-1}} \rio. 
 \end{eqnarray*}
By Lemma~\ref{convex}, 
 $$
2 \frac{F(i)}{\a \pi_i} - \frac{F(i-1)}{\a \pi_{i-1}} \le  \frac{F(i+1)}{\a
\pi_{i+1}}.
 $$
It follows that
 $$
\D^2 g_{f_{i4}} (i)  \le \D^2 g_{f_{1}} (1) +  2 \pi_{i+1} \leo  \frac{\bar{F}
(i+2)}{\a \pi_{i+1}} +  \frac{F(i+1)}{\a \pi_{i+1}} \rio = \D^2 g_{f_{1}} (1)+2 / \a
,
 $$
and the proof is completed by applying Lemma~\ref{lemdelta2g1}. \qed

\noindent{\bf Acknowledgements:} {We wish to thank a referee and an associate editor for extremely helpful
suggestions about the presentation of the paper. The research is} supported by the ARC Centre of 
Excellence for Mathematics and Statistics of Complex Systems.

 \end{document}